\documentclass[12pt,twoside,a4paper]{article}

\usepackage{german}
\usepackage{latexsym}
\usepackage{amsfonts}
\usepackage{amssymb}

\renewcommand{\baselinestretch}{0.55}
\input xy
\xyoption{all}

\def\dfrac{\displaystyle\frac}

\begin{document}
\vspace{7cm}
\title{Birational models for varieties of Poncelet curves}

\author{Matei Toma\\
\\
{\footnotesize Universit"at Osnabr"uck} \\
{\footnotesize Fachbereich Mathematik/Informatik}\\
{\footnotesize 49069~Os\-na\-br"uck, Germany}\\
{\footnotesize and}\\
{\footnotesize Institute of Mathematics of the Roumanian Academy} \\
{\footnotesize P.O.Box 1-764, RO-70700 Bucharest, Romania.}
}
\date{}
\maketitle
\thispagestyle{empty}
{\footnotesize
We consider Poncelet pairs $(S,C)$, where $S$ is a smooth conic
and $C$ is a degree$-c$ plane curve having the Poncelet property with respect
to $S$. We prove that for $c>4$ the projection $(S,C)\mapsto C$ is
generically one-to-one and use this to describe a birational model of the
variety of Poncelet curves for $c$ odd.
}

\section*{Introduction}
A pair $(S,C)$ consisting of a smooth conic $S$ and a curve $C$ of degree $c$
in the complex projective plane will be called a Poncelet pair if there exist
$c+1$ tangents to $S$ whose ${c+1\choose 2}$ intersection points lie on $C$.
$C$
will then be called Poncelet related to $S$ or simply a Poncelet curve.
A theorem of
Darboux (\cite{D}) states that in this case there is an infinity of such sets
of
$c+1$ tangents.

The Poncelet curves reappeared in a natural way in the study of stable vector
bundles on projective spaces, e.g. in \cite{B}, \cite{H}, \cite{BT}, \cite{NT}.
 Motivated by these developments, Trautmann gave a modern description of their
geometry in \cite{T}.

In this paper we consider the variety $Pon_c$ of Poncelet curves of degree $c$.
It is the image of the space of Poncelet pairs $(S,C)$ through the second
projection. We show that for $c\ge 5$, this projection is birational. This
fails for $c\le 3$, simply because  of dimension reasons, and has been proved
for $c=4$ by Le Potier via the study of divisors on the corresponding moduli
space of $rank -2$ semi-stable sheaves on $\Bbb P^2$, (\cite{L}).
As a corollary,
we get the rationality of $Pon_c$ for every $c$.

The generic injectivity of the above mentioned projection brings some evidence
for (and maybe could help proving) its generalizations to:
\begin{itemize}
\item Barth's morphism which associates to a stable $rank -2$  vector bundle of
type
$(0,c)$ on $\Bbb P^2$ its curve of jumping lines (\cite{B}), and
\item the restriction of the previous to Hulsbergen bundles. (The image will be
the variety of Darboux curves, cf. \cite{B}, \cite{ES}).
\end{itemize}
In the last paragraph we describe a projective birational model of $Pon_c$ when
$c$ is odd.

\section{Preliminaries}
For the proofs of the facts stated in this paragraph as well as for further
descriptions we refer the reader to \cite{T}.

Let $W$ be a 3-dimensional complex vector space, $S$ a smooth conic in $\Bbb
P^2=\Bbb P(W)$, $c$ a positive integer and $\Lambda$ a pencil of degree-$(c+1)$
divisors on $S$. We consider the intersection points of the tangents to $S$ at
the points of some member $D$ of
$\Lambda$. When $D$ moves in $\Lambda$, these intersection points describe a
curve $C=C(\Lambda)$ in $\Bbb P^2$. It is clear that the tangents at the base
points of $\Lambda$ will be components of $C$. We require that their
multiplicity in $C$ be equal to that of the base points in $D$. Then $C$ has
degree $c$, is called {\bf Poncelet related to} $S$ or simply a {\bf Poncelet
curve} and $(S,C)$ is called a {\bf Poncelet pair}.

That this definition is equivalent to that of the introduction follows from the

{\bf Theorem.}
{\it (Darboux) Let $S$ be a smooth conic and $C$ a curve of degree
$c$ in $\Bbb P^2$. If there exist $c+1$ tangents to $S$ such that their
intersection points lie on $C$, then $C$ is Poncelet related to $S$.
}

The base points of $\Lambda$ correspond to tangents to $S$ which are components
of $C(\Lambda)$ and the residual curve (obtained by eliminating these
components) is again Poncelet related to $S$ and associated to the pencil
obtained from $\Lambda$ by subtracting its base points.

For $\Lambda$ generic $C(\Lambda)$ is smooth.

Let's denote by $G_S(1,c+1)$ the Grassmannian of pencils of degree-$(c+1)$
effective divisors on $S$. The map
\begin{eqnarray}
G_S(1,c+1) &\to & \Bbb P(S^cW^*)=\Bbb P^{{c+2 \choose 2}-1}\nonumber\\
\Lambda &\mapsto & C(\Lambda)\nonumber
\end{eqnarray}
coincides with the Pl"ucker embedding (for a suitable choice of coordinates,
\cite{T}). Thus, there is a 1:1 correspondence between $G_S(1,c+1)$ and the
variety $Pon_{c,S}$ of curves $C$ which are Poncelet related to $S$.

Let $U$ be the open set of $\Bbb P(S^2W^*)$ which parameterizes smooth conics
in
$\Bbb P^2$, $I_U\subset U\times\Bbb P^2$ the tautological conic bundle
over $U$, $I^{(c+1)}_U$ the Hilbert scheme of degree-$(c+1)$ effective
divisors in the fibers of $I_U\to U$ (which is a $\Bbb P^{c+1}$-fibration over
$U$) and $G_U(1,c+1)$ the relative Grassmannian of lines in the fibers of
$I_U^{(c+1)}\to U$. We shall identify $G_U(1,c+1)$ with the subvariety
$Pon_{c,U}$ of Poncelet pairs in $U\times \Bbb P(S^cW^*)$. (The equations of
$Pon_{c,U}$ are given in \cite{T}).

We shall examine the second projection $pr_2: Pon_{c,U}\to\Bbb P(S^cW^*)$. Its
image is the variety of degree-$c$ Poncelet curves $Pon_c$.

We end this paragraph with some elementary remarks on Poncelet curves to be
used later:

{\it Remarks.} Let $S$ be a smooth conic.
\begin{enumerate}
\item[1.1] If $C$ is Poncelet related to $S$, then the only components of $C$
which may be multiple are tangent lines to $S$. (Apply Bertini's theorem to the
associated pencil, $\Lambda$.)
\item[1.2] If $C_1+C_2$ and $C_1+C_2'$ are two elements of $Pon_{c,S}$ such
that $C_2$ and $C_2'$ have no common component, then $C_1$ is Poncelet related
to $S$. (Start from a point on $C_1$, draw tangents to $S$ and continue this
procedure from the intersection points with $C_1$. The process stops after a
finite number of steps. If there exist intersection points of the drawn
tangents not on $C_1$, then these points lie on a common component of $C_2$ and
$C_2'$).
\item[1.3] Let $C_1, C_1+C_2$ be Poncelet related to $S$ with associated
pencils $\Lambda_1$ and $\Lambda$, respectively. If $\Lambda$ has no base
points
,
then the induced morphism $\varphi_\Lambda: S\to \Lambda^*\cong\Bbb P^1$
factorizes through $\varphi_{\Lambda'}:S\to\Lambda^{'*}$. In particular, $\deg
C_1+1$ divides $\deg(C_1+C_2)+1$.
\item[1.4] If $C$ is a singular conic Poncelet related to $S$, then one of its
components is tangent to $S$. (Follows from 1.3).
\end{enumerate}

\section{Generic injectivity of $pr_2$}
{\bf Theorem.}
{\it  When $c\ge 5$, the projection $Pon_{c,U}\to Pon_c$ is
birational.}

{\it Remarks.}
\begin{enumerate}
\item[2.1] The corresponding statement for $c=4$ (and more) is proven in
\cite{L}.
\item[2.2] $Pon_{c,U}$ is rational (\cite{HN}, Prop.2.2), hence so is
$Pon_c$ too.
\end{enumerate}
The theorem follows by a dimension count out of the following

{\bf Proposition.}
{\it Let $S,S'$ be distinct smooth conics. Then
\begin{enumerate}
\item[(A)] for $c\ge 3$, $\dim (Pon_{c,S}\cap Pon_{c,S'})\le c$.
\item[(B)] $Pon_{5,S}\cap Pon_{5,S'}$ is at most 4-dimensional at points
represented by smooth quintics.
\end{enumerate}
}

Indeed, let us estimate the dimension of
\linebreak
$\tilde{G}_S:=\overline{pr_2^{-1} (Pon_{c,S})\backslash
G_S(1,c+1)}$.  By part (A) of the Proposition,
$\dim \tilde{G}_S\le c+5$. Thus, for $c\ge 6$ we have $\dim pr_2(\tilde{G}_S)
\linebreak
< 2c=\dim Pon_{c,S}$ proving the theorem in this case.

When $c=5$, we use part (B) for a slightly modified $\tilde{G}_S$,
\linebreak
$\tilde{G}_S:=\overline{pr_2^{-1}(Reg_{c,S})\backslash G_S(1,c+1)}$, where
$Reg_{c,S}$ is the set of
\linebreak
smooth Poncelet curves related to $S$.

In order to prove the Proposition we look at the linear projection of
$Pon_{c,S}\cap Pon_{c,S'}$ with center $L_d:=\{\mbox{degree-$c$}$ curves in
$\Bbb P(W)$ allowing $d$ as a component$\} \subset \Bbb P(S^cW^*)$, for some
common tangent line $d$ to $S$ and $S'$.

Let $d$ be the tangent at $P$ to the smooth conic $S$. Then $Pon_{c,S}\cap L_d
\cong Pon_{c-1,S}$ corresponds to degree-$c$ pencils on $S$ having $P$ as base
point (or equivalently, to the Schubert variety on $G_S(1,c+1)$ of lines
contained in a hyperplane).

Let further $D$ be an element of $|{\cal O}_S(c+1)|$ containing $P$
and $L_D\subset
G_S(1,c+1)$ the Schubert variety of lines through $D$. Then $L_D\cong \Bbb P^c$
is linear in $\Bbb P(S^cW^*)$ and cuts on $L_d$ a $(c-1)$-dimensional linear
subspace, (the $L_{D-P}$ of $G_S(1,c)$). Thus $L_D\cup L_d$ spans a
$c$-dimensional linear subspace of $\Bbb P(S^cW^*)$.

Conversely, for every $c$-codimensional linear subspace $F$ which contains
$L_d$, these exists some $D$ as above such that $F\cap Pon_{c,S}=L_D\cup
(Pon_{c,S}\cap L_d)$.

For two distinct divisors $D,D'$ as above, $L_D$ and $L_{D'}$ cut themselves in
exactly one point lying on $L_d$. In particular, the projection of center $L_d$
displays $Pon_{c,S}\backslash L_d$ as a vector bundle over $|{\cal O}_S(c)|$.

Let now $V$ be an irreducible component of $Pon_{c,S}\cap Pon_{c,S'}$
not contained in $L_d$.

Let $\pi:V\backslash L_d\to|{\cal O}_S(c)|$ be the restriction to $V$
of the above
projection, $t:=\dim\pi(V\backslash L_d)$, and $f$ the dimension of a
generic fiber of $\pi$. Note that the fibers of $\pi$ are linear.

{\bf Lemma 2.1}
{\it When $c\ge 2$, up to finitely many exceptions, the fibers of
$\pi$ are at most $(c-1)$-dimensional.}

{\it Proof.} We choose a general hyperplane in $|{\cal O}_S(c)|$ and show that
the
fibers of $\pi$ over its points cannot be $c$-dimensional. Indeed, take $P_1$
on $S$ such that it doesn't belong to a common tangent to $S$ and $S'$.
The divisors containing $P_1$ form a hyperplane in $|{\cal O}_S(c)|$.
Take any such
divisor and add $P$ to it in order to obtain a divisor $D=P+P_1+P_2+\ldots
+P_c$ in $|O_S(c+1)|$.
Take further a pencil with base points $P_1,P_2,\ldots,P_{c-1}$ and containing
$D$. Its associated Poncelet curve is the union of the tangents at $P_1,\ldots,
P_{c-1}$ to $S$ and a line through the intersection point of $d$ with the
tangent to $S$ at $P_c$. By Remarks 1.2 and 1.3 not all such curves are
Poncelet related to $S'$.
\hfill$\Box$

\vspace{2ex}
{\it Proof of part (A) of the Proposition}

We argue by induction on $c$.

Let $c=3$ and suppose that $t+f>3$. Then $t\ge 2$ by Lemma 2.1. We pick then a
divisor $D=P+P_1+P_2+P_3$ on $S$ and look at the fiber of $\pi$ contained in
$L_D$.

We get a divisor $D'=P'+P'_1+P'_2+P'_3$ on $S'$ such that this fiber is exactly
$L_D\cap L_{D'}\backslash L_d$. ($L_{D'}$ is to be considered with respect to
$S'$).
Note that $P_1,P_2$ or $P_1,P'_2$ may be freely chosen so that this fiber be
non empty.

Denote by $d_i,d'_i$ the tangents at $P_i$ to $S$ and at $P'_i$ to $S'$,
respectively, $i\in\{1,2,3\}$.
We index our $P_i'-s$ so that $d\cap d_i=d\cap d'_i$.

Choose now $P_1,P_2'$ such that the point $d_1\cup d_2'$ does not lie on a
common tangent to $S$ and $S'$. In particular, it will not lie on both $d_3$
and $d'_3$, say not on $d'_3$.

The Poncelet curves in $L_{D'}\cap L_D$ touch $d'_2$ at its intersection points
with $d, d'_1$ and $d_3'$. (If one of the $P'_i$ had multiplicity $(k+1)$ in
$D'$ then the condition on the corresponding Poncelet curves would be to touch
$S$ at $P'_i$ with multiplicity $k$, by \cite{T}, Prop. 2.6). If we impose to
such a curve that $d_1$ be a component (which we may do since $f\ge 1$), we get
 that $d'_2$ must be a component too. By Remark 1.4 the third component of our
cubic is tangent to both $S$ and $S'$, so $f=1$ and $t=3$. This allows a free
choice on $P_3'$ too and we get that the tangents through $P_2'$ and $P_3'$ to
$S'$ form a conic which is Poncelet related to $S$. But this contradicts
Remark 1.4.

Let now $c>3$.

{\bf Case (A1)} $t>\dfrac{c}{2}$, $f\ge 1$.\\
We choose as above $D=P+P_1+\ldots+P_c$ and get in the same way a
$D'=P'+P'_1+\ldots+P'_c$ such that
$L_D\cap L_{D'} \backslash L_s\neq \emptyset$.
In doing this $P_1,\ldots, P_t$ are free.

We claim that the choice may be done such that $d_1$ intersects
$d'_2,\ldots,d'_t$ away from the vertices of the polygon $d,d'_1,\ldots,d'_c$
(with the obvious extension of notation).
If this were  not true then there would exist a point, say $P'_{t+1}$, in a
``bad position'', e.g. such that $d_1\cap d'_2\cap d'_{t+1}\neq \emptyset$.
We look
then at $P_2$ instead of $P_1$ and may find a second point $P'_{t+2}$, in a bad
position, and continue in this way. The assumption $t>\dfrac{c}{2}$ shows now
that our claim holds.

Next, from $f\ge 1$ it follows that we may degenerate a curve $C$ in our fiber,
so that $d_1$ and hence also $d'_2,\ldots,d'_t$ become its components. By
Remarks 1.2-4 the other components of $C$ are not all tangent lines to $S$,
otherwise we move the lines $d'_i$, $2\le i\le t$. In particular $t<c$.

If $f\ge 2$, then we may still move our curve $C$ (in its class), fact which
combined with Remark 1.3 shows the existence of a pair of degrees $(a,b)$ such
that $t-1\le a<b<c$ and $a+1$ divides $b+1$. But this would imply
$t\le\dfrac{c}{2}$, which proves that $f=1$ and $f+t\le c$.

{\bf Case (A2)} $f>\dfrac{c}{2}$\\
The proof is similar to that of Lemma 1.

We may first assume that the general $D=P+P_1+\ldots+P_c$ in the image of $\pi$
doesn't contain any tangency point of a common tangent to $S$ and $S'$
excepting $P$. Otherwise  we restrict our attention to the subvariety of
Poncelet curves containing this tangent and apply the induction hypothesis.

In the fiber of $\pi$ over $D$ we require that the tangents
$d_1,\ldots,d_{f-1}$ at $P_1,\ldots,P_{f-1}$ to $S$ be components of the
Poncelet curves $C$. Then by Remarks 1.3 and 1.4 there exists a pair of degrees
$(a,b)$ with $f-1\le a<b\le c$ and $a+1$
divides $b+1$. This forces $f=\dfrac{c+1}{2}$ proving part (A) of the
Proposition.
\hfill$\Box$

{\it Proof of part (B) of the Proposition}\\
Keeping our previous notations we assume that $t+f=c=5$.

{\bf Case (B1)}, $t=5$.\\
Since in this case $\pi$ is dominant we find for any choice of points
$P_1,\ldots,P_c$ on $S$, an element in $L_D\cap V$, where $D=P+P_1+\ldots+P_c$.
Take $P_1$ on $S$ and $P'_1$ the corresponding point on $S'$ (draw the tangent
$d'_1$ from $d\cap d_1$ to $S'$). Choose further a point $Q$ on $d'_1$ and draw
tangents from $Q$ to $S$. Let these tangents be $d_2,d_3$. A curve in $L_D\cap
V$ has to contain $Q$ and the intersection points of $d'_1$ with
$d,d'_2,\ldots, d'_c$. Thus this curve has $d'_1$ as a component. In the same
way, choosing $d_4$ through $d'_2\cap d_1$ fixes $d'_2$ as a component of $C$.
We may still move our point $P_5$, and since $d'_1\cup d'_2$ is not in
$Pon_{2,S}$, we get by Remarks 1.2 and 1.3 degrees $a,b$ with $3\le a<b\le 6$
and $a+1\mid b+1$. This is absurd so case (B1) cannot occur.

{\bf Case (B2)}, $t=4$.\\
As in case (A1) we  may fix the components $d_1,d'_2,\ldots,d'_t$ as we wish.
The last component of our Poncelet curve in not tangent to $S$ as already
remarked in (A1). Since we may consider $d'_2, d'_3$ fixed and move $d'_4$ we
get again a contradiction of Remarks 1.2 and 1.4. Case (B2) is thus excluded.

{\bf Case (B3)}, $t=3$, reduces itself to (A1).\\
{\bf Case (B4)} $t=2$.

We pick as usual two arbitrary tangent lines to $S$ and look at the curves in
$V$ containing them. They form an at least 1-dimensional subvariety. Since
these two tangents do not form a Poncelet conic with respect to $S'$ we get
degrees $3\le a<b\le c$ such that $a+1$ divides $b+1$. This implies $c\ge 7$
and excludes this case too.

{\bf Case (B5)}, t=1.\\
If we find some $D=P+P_1+\ldots+P_c$ in the ``image'' of $\pi$ containing two
points which don't lie on common tangents to $S$ and $S'$ we shall argue as in
(B4).

Let's assume then that $P_2,\ldots,P_c$ are fixed, all lying on common
tangents.
We have free choice on $P_1$. Imposing that $d_1$ become a component, we are
left with an $(f-1)$-dimensional family of curves $C$ in $L_{D-P_1}$, Poncelet
related to $S$, and such that $C+d_1$ belongs to $Pon_{c,S'}$.
Applying to this family the same
procedure with a different $P_1$, we get a 1-dimensional family of curves $C$
with the above
properties with respect to both choices of $P_1$. The curves $C$ do not
consist only of tangents to $S'$, since these would have to touch $S'$ only at
$P_2',\ldots,P_c'$. We thus get a contradiction of Remarks 1.2-4 when applied
to $Pon_{c,S'}$.

{\bf Case (B6)}, $t=0$.\\
In this case there exist a divisor $D=P+P_1+\ldots+P_c$ on $S$ such that
$V=L_D$. The argument from (B4) shows that
at most one point, $P_1$ say, is not on a common tangent to $S$ and $S'$. By
requiring that $d_3,\ldots, d_6$ be components we reduce ourselves to the case
$c=2=f$, $D=P+P_1+P_2$. But this is excluded by showing that $d_2$ has to be a
component of $C$.

If all points of $D$ are tangency points of common tangents to $S$ and $S'$,
then either two of them are multiple in $D$ or one of them has multiplicity
bigger than 3 in $D$. In both cases the associated Poncelet curves are all
singular by \cite{T}, Proposition 5.1.
\hfill$\Box$

\section{A projective birational model of $Pon_c$ for $c$ odd}
For odd $c$, $G_U(1,c+1)$ is the relative Grassmannian of 2-dimen\-sional
subspaces of a rank-$(c+2)$ vector bundle over $U$. This vector bundle and thus
also $G_U(1,c+1)$ may be extended to the whole of $\Bbb P(S^2W^*)$.
When $c\ge 5$ one obtains a birational map from this extended relative
Grassmannian to $\overline{Pon_c}$.

In this paragraph we describe the linear system which induces this birational
map and its base locus as a set. (Note that the knowledge of the full scheme -
structure of this base-locus would allow one to compute the degree of
$\overline{Pon_c}$). Our method is to consider the stable rank-two vector
bundles on $\Bbb P(W)$ associated to Poncelet curves and their curves of
jumping lines.

Throughout this paragraph $c$ is assumed to be odd and bigger than 4.

The extension of $G_U(1,c+1)$ to $\Bbb P(S^2 W^*)$ is obvious. Take $I\subset
\Bbb P(S^2W^*)\times \Bbb P(W)$ the incidence variety ``points of conics'',
$p_1:\Bbb P(S^2W^*)\times \Bbb P(W)\to \Bbb P(S^2W^*)$ the first projection,
${\cal V}:=p_{1,*}{\cal O}_I\left(0,\dfrac{c+1}{2}\right)$ and $G$ the relative
Grassmannian of 2-dimen\-sional subspaces in the fibers of ${\cal V}$. $G$ is
then a compactification of $G_U(1,c+1)$.

Consider now a Poncelet pair $(S,C)$ with $S$ smooth and $C=C(\Lambda)$ with
base-point-free pencil $\Lambda\subset\left|{\cal O}_S\left(\dfrac{c+1}{2}
\right)\right|$. (In this paragraph we use the notation
${\cal O}_S\left(\dfrac{c+1}{2}\right):={\cal O}_{\Bbb
P(W)}\left(\dfrac{c+1}{2}\right)\mid_S$). $\Lambda$ induces a surjective
morphism
$$
{\cal O}_{\Bbb P(W)}^2\to {\cal O}_S\left(\dfrac{c+1}{2}\right).
$$
Let $F$ be its kernel. Then $F$ is a stable rank-2 vector bundle on $\Bbb P(W)$
with $c_1(F)=-2$, $c_2(F)=c+1$. and its jumping lines are exactly the lines
joining the points of some divisor of $\Lambda$ (see \cite{T}, 4.1). Thus, the
curve $C'\subset \Bbb P(W^*)$ of jumping lines of $F$ is Poncelet related to
the conic $S'$ dual to $S$.

We shall make $F$ fit into a flat family of coherent rank-2 sheaves over $G$,
and examine the birational map associating the curves of jumping lines to the
fibers of this family.

We start with some  Lemmata.

{\bf Lemma 3.1.}
{\it Let $0\to E'\to {\cal O}^N_{\Bbb P^2} \to E''\to 0$ be an
exact sequence with $E',E''$ coherent torsion-free, indecomposable sheaves on
$\Bbb P^2$ and $c_1(E')=-1$. Then $N\le 3$ and $E',E''$ are slope-stable.}

{\it Proof.}
$
E''$ torsion-free implies $E'$ locally-free.
If $E'$ were not slope-stable, there would exist a locally free subsheaf $F$ of
$E'$ with $c_1(F)=0$. Standard arguments on slope-stable vector bundles (cf.
\cite{K}, V 8.3) show that $F$ would be a direct summand in $E'$,
contradicting the
hypothesis.

Similarly, $E''$ must be slope-stable.

Let $r':=rank\; E'$, $r'':=rank\; E$.

Suppose $r'=1$. Then $E'={\cal O}_{\Bbb P^2}(-1)$ and the morphism ${\cal
O}_{\Bbb P^2}(-1)
\linebreak
\to  {\cal O}^N_{\Bbb P^2}$ is induced by $N$ sections in
${\cal O}_{\Bbb P^2} (1)$. Since $\Gamma({\cal O}_{\Bbb P^2}(1))$ is
3-dimensional we may choose when $N>3$ a basis for ${\cal O}_{\Bbb P^2}(1)^N$
such that $N-3$ components of ${\cal O}_{\Bbb P^2}\to {\cal O}_{\Bbb
P^2}(1)^N$ vanish. But then $E''$ would split.

The case $r''=1$ is treated in the same way.

We are left with the situation $r'>1$, $r''>1$. But now the Bogomolov
inequality gives
$$
\renewcommand{\arraystretch}{2.5}
\begin{array}{l}
c_2(E')\ge \dfrac{r'-1}{2r'}\cdot c_1(E')^2>0\\
\\
c_2(E'')\ge \dfrac{r''-1}{2r''}\cdot c_1(E'')^2>0
\end{array}
$$
hence $c_2({\cal O}^N_{\Bbb P^2})=-1+c_2(E')+c_2(E'')>0$, a contradiction.
\hfill$\Box$

{\bf Lemma 3.2.}
{\it Let ${\cal K}$ be the kernel of the natural surjective morphism
$$
p_1^*{\cal V}\to {\cal O}_I\left(0,\dfrac{c+1}{2}\right)
$$
on $\Bbb P(S^2W^*)\times\Bbb P(W)$. Then ${\cal K}$ is locally free and
its restrictions to the fibers of $p_1$ are slope-stable.}

{\it Proof.} The restriction to a fiber of $p_1$ over some point $s$
representing a conic $S$ gives an exact sequence
$$
0\to{\cal K}\mid_{p^{-1}(s)}\to\Gamma\left({\cal O}_S
\left(\dfrac{c+1}{2}\right) \right)
\otimes {\cal O}_{\Bbb P(W)}\to {\cal O}_S\left(\dfrac{c+1}{2}\right)\to 0.
$$
It is enough to check that $K:={\cal K}\mid_{p^{-1}(s)}$ is locally free and
slope stable.

Checking the locally freeness is easy (compare $K$ to its double dual or
just use the fact that ${\cal K}$ appears from an elementary transformation!).

Suppose that $K$ were not slope-stable.

Since $c_1(K)=-2$, there would exist a destabilizing subsheaf $K'$ of $K$ with
$c_1(K')=-1$.
By Lemma 3.1. $K'$ would have at most rank 2, hence the rank of $K$ could not
exceed 4.

But rank $K=c+2$ and the Lemma follows from our assumption on $c$.
\hfill$\Box$

Let ${\cal S}$ and ${\cal Q}$ be the tautological sub- and quotient bundle
on the relative Grassmannian $G$ of 2-dimensional subspaces in the fibers of
${\cal V}$. We denote by ${\cal S}_{\Bbb P(W)}$, ${\cal Q}_{\Bbb P(W)}$,
${\cal K}_G$, ${\cal V}_{G\times\Bbb P(W)}$ the pullbacks to $G\times \Bbb
P(W)$ of ${\cal S, Q, K}$ and ${\cal V}$, respectively. (Subscript will be used
in the sequel to indicate pullbacks through obvious maps).
Denote further by $\alpha$ the composite morphism
$$
{\cal Q}_{\Bbb P(W)}^{\vee}\to {\cal V}_{G\times\Bbb P(W)}^{\vee}\to
{\cal K}_G^{\vee}
$$
on $G\times\Bbb P(W)$ and by ${\cal E}$ its cokernel.

{\it Remark.} One sees immediately that the fiber of ${\cal E}$ over a point of
$G$ represented by a Poncelet pair $(S,C(\Lambda))$ with $S$ smooth and
$\Lambda$ base-point-free is just the dual of the kernel $F$ of the natural
morphism
$$
{\cal O}_{\Bbb P^2(W)}^2\to{\cal O}_S\left(\dfrac{c+1}{2}\right).
$$
In particular, it is stable and locally free.

We claim that ${\cal E}$ is flat over $G$. This is a consequence of a local
flatness criterion (\cite{M}; 2.2.5) and the following

{\bf Lemma 3.3}
{\it  The restrictions of $\alpha$ to the fibers of $G\times\Bbb
P(W)\to G$ are injective and their cokernels are stable as soon as they are
torsion-free.}

{\it Proof.} Fix a conic $S$ and a $c$-dimensional subspace $A$ of
$\Gamma
\linebreak
\left({\cal O}_S \left(\frac{c+1}{2}\right)\right)^* \cong \Gamma(\Bbb P(W),
 K^{\vee})$.  The slope-stability of $K$ and Lem\-ma 3.1 imply that
the natural morphism $s:A\otimes{\cal O}_{\Bbb P(W)}\to K^{\vee}$ is injective.

If Coker $s$ is torsion-free but not stable, it will admit a rank-1 subsheaf
$E'$ with torsion-free quotient and $c_1(E')=1$, $c_2(E')\le\frac{c-1}{2}$. In
particular, $\dim\: Ext^1(E',
{\cal O}_{\Bbb P(W)})\le\frac{c-1}{2}$. We get a
subsheaf $K'$ of $K^{\vee}$ with $c_1(K')=1$, $rank(K')=c+1$,
and which sits in an
exact sequence
$$
0\to A\otimes{\cal O}_{\Bbb P(W)}\to K'\to E'\to 0.
$$

This leads to a contradiction of the stability of $K$ in view of the following
simple fact of homological algebra:

{\bf Lemma 3.4}
{\it Let $A$ be a finite-dimensional $\Bbb C$-vector space,
${\cal C}$ a coherent sheaf on a compact complex manifold $X$ and
$$
0\to A^*\otimes{\cal O}_X\to {\cal B}\to {\cal C}\to 0
$$
an extension given by an element
$$
\eta\in Ext^1({\cal C}, A^*\otimes{\cal O}_X)\cong Hom\:(A, \: Ext^1({\cal
C},{\cal O}_X))
$$
such that $A=A'\oplus A''$ and $\eta\mid_{A''}=0$. Then
$$
\eta\mid_{A'}\in Hom\:(A', \: Ext^1({\cal C},{\cal O}_X))\cong Ext^1({\cal
C},A^{'*}\otimes{\cal O}_X)
$$
gives rise to a commutative diagram

$$
\xymatrix{
0 \ar[r] & A^{'*}\otimes {\cal O}_X \ar[d] \ar[r]
& {\cal B}'\ar[d] \ar[r]
& {\cal C} \ar[d] \ar[r]
& 0\\
0 \ar[r]
& A^*\otimes {\cal O}_X \ar[r]
& {\cal B} \ar[r]
& {\cal C} \ar[r]
& 0 }
$$

with injective vertical maps.}
\hfill$\Box$

Associating to a fiber of ${\cal E}$ over $G$ its divisor of jumping lines
defines a rational map
$G \to \Bbb P(S^cW) $ in the usual way. Indeed, we
consider the standard diagram

$$
\xymatrix{
\Bbb F \ar[d] \ar[r] & \Bbb P(W^*)\\
\Bbb P(W) & }
$$

where ${\Bbb F}$ is the incidence variety ``points on lines'', then take its
product with $G$

$$
\xymatrix{
G\times \Bbb F \ar[d]_p \ar[r]^q & G\times \Bbb P(W^*)\\
G\times \Bbb P(W) & }
$$

and look at the relative theta characteristic associated to ${\cal E}$,
$$
\Theta ({\cal E}):= R^1q_* p^*({\cal E}(-2)).
$$

One sees immediately that, for the torsion free fibers, ${\cal E}_g$ of ${\cal
E}$, the corresponding fibers of $\Theta({\cal E})$ are the usual
theta-characteristics associated to a stable sheaf on $\Bbb P^2$ with even
first Chern
class (cf. \cite{B}).
In particular, their supports are the curves of jumping lines
of the ${\cal E}_g-s$. Let $B$ be the subset of $G$ over which ${\cal E}$ has
nontrivial torsion in its fibers. Then we get a morphism
$$
G\backslash B\to \Bbb P(S^cW),
$$
(see \cite{Ma}).

 From the defining sequence of ${\cal E}$ one obtains the following resolution
of $\Theta({\cal E})$ on $G\times \Bbb P(W^*)$
$$
0\to{\cal Q}^{\vee}_{\Bbb P(W^*)}(-1)\stackrel{\varphi}{\to}\left( R^{1}
p_{1*}{\cal K}^{\vee}\right)_{G\times\Bbb P(W^*)}\to \Theta({\cal E})\to 0.
$$

In order to obtain the ``support'' $Z$ of $\Theta({\cal E})$ one takes the
determinant of $\varphi$ and twists it correspondingly to get:
$$
0\to{\cal L}^{-1}\boxtimes {\cal O}_{\Bbb P(W^*)}(-c)\to{\cal O}_{G\times\Bbb
P(W^*)}\to{\cal O}_Z\to 0
$$
where
$$
{\cal L}:=\det {\cal Q}\otimes {\cal O}_{\Bbb P(S^2W^*)}
\left(-\dfrac{(c-1)(c-3)}{8}\right)_G.
$$
Notice that $Z$ is not flat over $G$. In fact, its fibers over $B$ are
2-dimensional. Thus the rational map $G \to \Bbb P(S^cW) $
above is given by a linear subsystem of $|{\cal L}|$ and has base locus $B$.
The
structural sheaf of this base locus may be recovered by taking the push-down
on $G$ of the previous exact sequence twisted by (-3). One obtains for it the
presentation
$$
{\cal L}^{-1}\otimes H^0({\cal O}_{\Bbb P(W^*)}(c))^*\to {\cal O}_G \to R^2
\pi_*({\cal O}_Z(-3))\to 0,
$$
where $\pi:G\times \Bbb P(W^*)\to G$ is the projection.

Finally we want to identify $B$.

Let $\psi:\Bbb P(W^*)\times \Bbb P(W^*)\to \Bbb P(S^2W^*)$ be the map induced
by the product, $p':\Bbb P(W^*)\times\Bbb P(W^*)\times \Bbb P(W)\to \Bbb
P(W^*)\times \Bbb P(W^*)$ the projection,
$$
I_i=\{(l_1,l_2,x)\in \Bbb P(W^*)\times \Bbb P(W^*)\times \Bbb P(W)\mid x\in
l_i\},\; i\in\{1,2\},
$$
the incidence varieties and ${\cal W}:=p_*'({\cal O}_{I_1}(0,-1,
\frac{c-1}{2}))$.

{\bf Proposition.}
{\it ${\cal W}$ is a subbundle of $\psi^*{\cal V}$ and the
natural map from the Grassmannian of 2-dimensional subspaces in the fibers of
${\cal W}$ to $G$ is an embedding whose image is $B$. In particular, $B$ is
smooth of dimension $c+5$.}

{\it Proof.} We examine when a fiber of ${\cal E}$ over $G$ admits nontrivial
torsion. Fix a conic $S$ and a 2-dimensional subspace $L$ of $\Gamma\left(
{\cal O}_S \left(\frac{c+1}{2}\right)\right)$. If $K$ and $E$ are the fibers
of ${\cal K}$ and ${\cal E}$ over the corresponding point of $G$ we have exact
sequences:
$$
\begin{array}{l}
0\to \left(\Gamma\left({\cal O}_S\left(\dfrac{c+1}{2}\right)\right)
\biggl/L\right)^* \otimes {\cal O}_{\Bbb P(W)}\to K^{\vee}\to E\to 0,\\
o\to E^{\vee}\to L\otimes {\cal O}_{\Bbb P(W)}\stackrel{\epsilon}{\to}{\cal
O}_S
\left(\dfrac{c+1}{2}\right)\to {\cal E}xt^1(E,{\cal O}_{\Bbb P(W)})\to 0.
\end{array}
$$
where ${\epsilon}$ is the evaluation morphism.
Since $K$ is slope-stable, $E$ may
only have purely 1-dimensional torsion, and this happens exactly when
${\cal E}xt^1(E,{\cal O}_{\Bbb P(W)})$
has 1-dimensional support or equivalently,
when ${\epsilon}$ fails to be ``generically surjective on some component of
$S$''. More precisely, if this is the case, then $S=S_1+S_2$ and the composite
morphism
$$
L\otimes {\cal O}_{\Bbb P(W)}\left(\dfrac{c+1}{2}\right)
\stackrel{\epsilon}{\to}
{\cal O}_S\left(\dfrac{c+1}{2}\right)\to{\cal O}_{S_2}
\left(\dfrac{c+1}{2}\right)
$$
vanishes. Then ${\epsilon}$ factors through ${\cal
O}_{S_1}\left(\dfrac{c+1}{2}\right)$ and the pullback of $L$ through $\psi$
will be a 2-dimensional subspace in a fiber of ${\cal W}$ over $\Bbb
P(W^*)\times \Bbb P(W^*)$.
Indeed, pushing the exact sequence
$$
0\to {\cal O}_{I_1}(0,-1,-1)\to\psi^*{\cal O}_I\to {\cal O}_{I_2}\to 0
$$
twisted by $\left(\frac{c+1}{2}\right)$ from $\Bbb P(W^*)\times\Bbb P(W^*)
\times\Bbb P(W)$ down to $\Bbb P(W^*)\times \Bbb P(W^*)$ gives
$$
0\to {\cal W}\to\psi^*{\cal V}\to p_*'({\cal
O}_{I_2}\left(0,0,\left(\frac{c+1}{2}\right)\right)\to 0.
$$

The claims of the Proposition are now easy to check; that the considered
morphism is an embedding follows e.g. from a computation of its Jacobian
matrix.

\hfill$\Box$

\markboth{}{}
\markright{}
{\footnotesize
\renewcommand{\baselinestretch}{0.3}

}
\end{document}